\documentclass[final,1p]{elsarticle}
\usepackage{graphics}
\usepackage{graphicx}
\usepackage{epsfig,bm}
\usepackage{amssymb,amsmath}

\journal{Nonlinear Analysis: Real World Applications - accepted.}

\newtheorem{definition}{Definition}
\newtheorem{prop}{Proposition}
\newtheorem{thm}{Theorem}
\newtheorem{cor}{Corollary}
\newtheorem{ex}{Example}
\newtheorem{lem}[thm]{Lemma}
\newdefinition{rem}{Remark}
\newproof{pf}{Proof}

\begin{document}

\begin{frontmatter}
\title{Non-existence of periodic solutions in fractional-order dynamical systems and a remarkable difference between integer and fractional-order derivatives of periodic functions \tnoteref{thanks}}
\tnotetext[thanks]{This work was supported by CNCSIS-UEFISCSU, project number PN-II-RU-PD-145/2010 (Advanced impulsive and fractional-order neural network models).}

\author[IeAT,UVT]{Eva Kaslik}
\ead{ekaslik@gmail.com}
\author[ERAU]{Seenith Sivasundaram\corref{cor}}
\ead{seenithi@gmail.com}
\cortext[cor]{Corresponding author}

\address[IeAT]{Institute e-Austria Timisoara, Bd. V. Parvan nr. 4, room 045B, 300223, Timisoara, Romania}

\address[UVT]{Dept. of Mathematics and Computer Science, West University of
Timisoara, Bd. V. Parvan nr. 4, 300223, Romania}

\address[ERAU]{Department of Mathematics, Embry-Riddle Aeronautical University, Daytona Beach, FL 32114, USA}


\begin{abstract}
Using the Mellin transform approach, it is shown that, in contrast with integer-order derivatives, the fractional-order derivative of a periodic function cannot be a function with the same period. The three most widely used definitions of fractional-order derivatives are taken into account, namely, the Caputo, Riemann-Liouville and Grunwald-Letnikov definitions. As a consequence, the non-existence of exact periodic solutions in a wide class of fractional-order dynamical systems is obtained. As an application, it is emphasized that the limit cycle observed in numerical simulations of a simple fractional-order neural network cannot be an exact periodic solution of the system.
\end{abstract}

\end{frontmatter}

\section{Introduction}

Although fractional calculus has a more than 300 year long history, its practical applications are just a recent focus of interest. In the past decade, scientists and engineers became aware of the fact that the description of some phenomena is more accurate when the fractional derivative is used. It has been recently found that many systems in interdisciplinary fields (such as acoustics, mechanics, electromagnetism, heat transfer, electrical circuits, signal processing, system identification, control and robotics, chemistry, biology, physics, economy and finance) can be successfully described by fractional differential equations.

Applications of fractional calculus and fractional-order differential equations include, but are not limited to modeling real world phenomena such as: dielectric relaxation phenomena in polymeric materials \cite{Reyes-Melo}, transport of passive tracers carried by fluid flow in a porous medium in groundwater hydrology \cite{Schumer_Benson}, viscoelastic behavior \cite{Heymans_Bauwens}, transport dynamics in systems governed by anomalous diffusion \cite{Henry_Wearne,Metzler}, self-similar processes such as protein dynamics \cite{Nonnenmacher}, long-time memory in financial time series \cite{Picozzi_West} using fractional Langevin equations \cite{Bashir_Langevin}, etc. In recent years, even fractional-order models of happiness \cite{Song_Xu_Yang} and love \cite{Gu_Xu} have been developed, and they are claimed to give a better representation than the integer-order dynamical systems approach.

Highly remarkable scientific books which provide the main theoretical tools for the qualitative analysis of fractional-order dynamical systems, and at the same time, show the interconnection as well as the contrast between classical differential equations and fractional differential equations, are \cite{Kilbas,Lak,Podlubny}.


The existence of periodic solutions is a often a desired property in dynamical systems, constituting one of the most important research directions in the theory of dynamical systems, with applications ranging from celestial mechanics to biology and finance. The existence of weighted pseudo-almost periodic solutions of fractional-order differential equations has been investigated in \cite{Agarwal_Andrade_Cuevas,Cao_Yang_Huang_2011,Debbouche_Borai}. Moreover, several results concerning the existence of solutions of periodic boundary problems for fractional differential equations have been recently reported in \cite{Bai_2011,Belmekki-Nieto,Nieto_1,Wei}.

In this paper, we will first show that, unlike in the case of the integer-order derivative, the fractional-order derivative of a periodic function cannot be a function with the same period. As a consequence, we obtain the non-existence of exact periodic solutions for a wide class of fractional-order differential systems. Finally, we exemplify the theoretical results in the framework of fractional-order neural network models.

\section{Preliminaries}

\subsection{Fractional order derivatives}

In general, three different definitions of fractional derivatives are widely used: the Grunwald-Letnikov derivative, the Riemann-Liouville derivative and the Caputo derivative. These three definitions are in general non-equivalent. However, the main advantage of the Caputo derivative is that it only requires initial conditions given in terms of integer-order derivatives, representing well-understood features of physical situations and thus making it more applicable to real world problems.

Let $C^n(I,\mathbb{R})$ denote the space of $n$-times continuously differentiable functions on the real interval $I$.

\begin{definition}
Let $\alpha\in(0,\infty)\setminus \mathbb{N}$ and the function $g\in C^n([a,b],\mathbb{R})$. We define:
\begin{itemize}
\item the Riemann-Liouville fractional-order derivative of order $\alpha$ of $g$, given by
$$
^{RL}D_{a+}^{\alpha}g(t)=\frac{1}{\Gamma(n-\alpha)}\left(\frac{d}{dt}\right)^n\int_
a^t(t-s)^{n-\alpha-1}g(s)ds,
$$
where $n=[\alpha]+1$.
\item the Caputo fractional-order derivative of order $\alpha$ of $g$, given by
$$
^CD_{a+}^{\alpha}g(t)=\frac{1}{\Gamma(n-\alpha)}\int_
a^t(t-s)^{n-\alpha-1}g^{(n)}(s)ds,
$$
where $n=[\alpha]+1$.
\item the Grunwald-Letnikov fractional-order derivative of order $\alpha$ of $g$, given by
$$
^{GL}D_{a+}^{\alpha}g(t)=\lim_{h\rightarrow 0}h^{-\alpha}\sum_{r=0}^{\left[(t-a)/h\right]}(-1)^r\left(\begin{array}{c}
                                                                                          \alpha \\
                                                                                          r
                                                                                        \end{array}
\right)g(t-rh)
$$
\end{itemize}
\end{definition}

For simplicity, when $a=0$, we will drop the subscript $"a+"$.

The following proposition (see \cite{Kilbas,Podlubny}) expresses the relationship between the three types of fractional derivatives.

\begin{prop}\label{prop.relation}
Let $\alpha\in(0,\infty)\setminus \mathbb{N}$ and $g:(0,\infty)\rightarrow$ a function of class $C^n$. Then
\begin{equation}\label{rel.Caputo.RL}
^{RL}D^\alpha g(t)=^{GL}D^\alpha g(t)=^{C}D^\alpha g(t) +\sum_{k=0}^{n-1}\frac{g^{(k)}(0^+)}{\Gamma(k-\alpha+1)}t^{k-\alpha}
\end{equation}
where $n=[\alpha]+1$.
\end{prop}

\subsection{The Mellin transform and its properties}

In the following, we will recall the definition of the Mellin transform and a few properties that will be helpful in the proof of our main result. We refer to \cite{Podlubny,Fikioris,Flajolet} and the references therein, for an overview of the Mellin transform and its applications.

\begin{definition}
The Mellin transform of a locally Lebesgue integrable function $g:[0,\infty)\rightarrow\mathbb{C}$ is defined by
$$\mathcal{M}(g)(z)=\int_0^\infty g(t) t^{z-1}dt.$$
The largest open vertical strip of the complex plane, of the form  $$S_g=\{z\in\mathbb{C}:~a<\Re(z)<b\},$$ in which the integral converges, is called the strip of analyticity (or fundamental strip) of the Mellin transform.
\end{definition}

\begin{prop}[Mellin transform of the convolution]\label{prop.Mellin.convolution}
Considering the Mellin convolution $g\ast h$ of two functions $g,h:[0,\infty)\rightarrow\mathbb{C}$ defined by
$$(g\ast h)(t)=\int_0^\infty g(ts)h(s)ds,$$
the following equality holds:
$$\mathcal{M}(g\ast h)(z)=\mathcal{M}(g)(z)\cdot\mathcal{M}(h)(1-z)$$
for any $z$ from the fundamental strip of $\mathcal{M}(g)$, such that $1-z$ belongs to the fundamental strip of $\mathcal{M}(h)$.
\end{prop}

\begin{prop}[Inversion of the Mellin transform]\label{prop.Mellin.inversion}
Let $g:[0,\infty)\rightarrow\mathbb{C}$ be an integrable function, such that its Mellin transform $\mathcal{M}(g)$ has the fundamental strip $S_g=\{z\in\mathbb{C}:~a<\Re(z)<b\}$. If $c\in (a,b)$ is such that $\mathcal{M}(g)(c+is)$ is integrable,  then the following equality holds:
$$\frac{1}{2i\pi}\int_{c-i\infty}^{c+i\infty}\mathcal{M}(g)(z)t^{-z}dz=g(t),$$
almost everywhere on $(0,\infty)$.
\end{prop}

\section{Main result}

For completeness, we first prove the following standard mathematical result:

\begin{lem}\label{lem.per}
Let $n\in\mathbb{N}$ and $T>0$. If $x:(0,\infty)\rightarrow\mathbb{R}$ is a non-constant $T$-periodic function of class $C^n$ on $(0,\infty)$, then for any $k\in\mathbb{N}$, $k\leq n$, the $k$-th order derivative $x^{(k)}$ is also a non-constant $T$-periodic function.
\end{lem}

\begin{pf}
Let $k\in\mathbb{N}$, $k\leq n$. If $x:(0,\infty)\rightarrow\mathbb{R}$ is a non-constant $T$-periodic function, it follows that
\begin{equation}\label{eq.per}
x(t+T)=x(t)\qquad\forall~t\in(0,\infty).
\end{equation}
Differentiating $k$ times in this equality, with respect to $t$, it follows that
\begin{equation}\label{eq.per.k}
x^{(k)}(t+T)=x^{(k)}(t)\qquad\forall~t\in(0,\infty),
\end{equation}
which means that $x^{(k)}$ is also $T$-periodic.

For the second part of the proof, we will proceed by \emph{reductio ad absurdum}.  Assuming that there exists $k\in\mathbb{N}$, $k\leq n$, such that $x^{(k)}(t)$ is a constant function, we obtain that $x(t)$ is a $k$-th degree polynomial function
$$x(t)=c_k t^k +c_{k-1}t^{k-1}+...+c_1t+c_0,$$
where $c_i\in\mathbb{R}$, for $i=\overline{0,k}$. Denoting $$y(t)=x(t)-c_0=c_k t^k +c_{k-1}t^{k-1}+...+c_1t,$$
it is clear from the fact that $x(t)$ is $T$-periodic, that the function $y(t)$ is also $T$-periodic. Moreover, we note that $y(0)=0$, and from $T$-periodicity, it follows that
$$y(mT)=y(0)=0\qquad\forall m\in\mathbb{N}.$$
Since $T>0$, this means that the $k$-th degree polynomial function $y(t)$ has an infinity of roots: $\{mT,~m\in\mathbb{N}\}$, which is absurd. Hence, our assumption is false, and therefore, for any $k\in\mathbb{N}$, $k\leq n$, the function $x^{(k)}(t)$ is non-constant. \qed
\end{pf}

\begin{rem}
In the proof of the previous Lemma, when we differentiate equality (\ref{eq.per}) with respect to $t$, we implicitly rely on the well-known chain rule for the integer-order derivative:
$$[f(g(t))]'=g'(t)\cdot f'(g(t)).$$
This allows us to differentiate the left hand side of (\ref{eq.per}), leading to:
$$\frac{d}{dt}x(t+T)=\frac{d}{dt}(t+T)\cdot x'(t+T)=1\cdot x'(t+T)=x'(t+T),$$
and hence:
$$x'(t+T)=x'(t)\qquad\forall t\in(0,\infty).$$
By successive differentiations, we obtain equality (\ref{eq.per.k}).

However, it is well known that there is no simple chain rule for fractional-order derivatives (see, for example, section 2.7.3 in \cite{Podlubny}). Hence, applying the fractional derivative operator $D^\alpha$, with $\alpha\in (0,\infty)\setminus\mathbb{N}$, to both sides of the equation
(\ref{eq.per}), will clearly not lead us to
$$(D^\alpha x)(t+T)=(D^\alpha x)(t)\qquad\forall~t\in(0,\infty).$$
\end{rem}

In the following, we present the main result of this paper.

\begin{thm}\label{thm.main}
Let $\alpha\in (0,\infty)\setminus\mathbb{N}$ and $n=[\alpha]+1$. If $x:(0,\infty)\rightarrow\mathbb{R}$ is a non-constant $T$-periodic function of class $C^n$, then its Caputo derivative $^CD^\alpha x$ cannot be a $T$-periodic function.
\end{thm}

\begin{pf}
Assume that there exists a non-constant periodic function $x(t)$ of period $T>0$, $n$-times continuously differentiable on $(0,\infty)$, such that its Caputo derivative is $T$-periodic as well, i.e.
$$^CD^\alpha x(t)=^CD^\alpha x(t+T)\qquad\forall~ t\geq 0.$$
Applying the definition of the Caputo derivative, we get:
$$\int_0^t(t-s)^{n-\alpha-1}x^{(n)}(s)ds=\int_0^{t+T}(t+T-s)^{n-\alpha-1}x^{(n)}(s)ds,$$
where $n=[\alpha]+1$.
Making the change of variables $s+T=s'$ in the integral from the left hand-side of the equality, and taking into account that $x^{(n)}$ is $T$-periodic as well (see Lemma \ref{lem.per}), we obtain:
$$\int_T^{t+T}(t+T-s)^{n-\alpha-1}x^{(n)}(s)ds=\int_0^{t+T}(t+T-s)^{n-\alpha-1}x^{(n)}(s)ds,$$
which is equivalent to
$$\int_{0}^{T} (t+T-s)^{n-\alpha-1}x^{(n)}(s)ds=0 \qquad,~\forall t> 0.$$
Next, we make the change of variables $s=T-ts'$ and obtain
$$\int_0^{\frac{T}{t}} (t+ts)^{n-\alpha-1}x^{(n)}(T-ts)tds=0\qquad,~\forall~t> 0$$
and, simplifying $t^{n-\alpha}$, this is equivalent to
$$\int_0^{\frac{T}{t}} (1+s)^{n-\alpha-1}x^{(n)}(T-ts)ds=0\qquad,~\forall~t> 0.$$
Defining the function
$$h(u)=\left\{
         \begin{array}{ll}
           x^{(n)}(T-u) & ,~\textrm{if } u\in [0,T] \\
           0 & ,~\textrm{if } u>T
         \end{array}
       \right.,
$$
the previous equation can be written as
$$\int_0^{\infty} (1+s)^{n-\alpha-1}h(ts)ds=0\qquad\forall~t> 0,$$
which translates into the fact that the Mellin convolution $h\ast g$ of the functions $h(t)$ and $g(t)=(1+t)^{n-\alpha-1}$ is equal to $0$.
Applying the Mellin transform to this equality, it follows that
$$\mathcal{M}(h\ast g)(z)=0,$$
for any $z$ where this Mellin transform is defined.

Using the formula for the Mellin transform of the convolution given by Proposition \ref{prop.Mellin.convolution}, we get
$$H(z)G(1-z)=0,$$
where $H=\mathcal{M}(h)$ and $G=\mathcal{M}(g)$. This equality is true for any $z\in\mathbb{C}$ such that $z\in S_h$ and $1-z\in S_g$, where $S_h$ and $S_g$ are the fundamental strips of analiticity of the Mellin transforms of $h$ and $g$.

It can be easily seen that the Mellin transform of $g$ is defined only on the strip $S_g=\{z\in\mathbb{C}:0<\Re(z)<\alpha-n+1\}$ and
$$G(z)=\int_0^\infty g(t)t^{z-1}dt=\int_0^\infty (1+t)^{n-\alpha-1}t^{z-1}dt=\frac{\Gamma(\alpha-n+1-z)\Gamma(z)}{\Gamma(\alpha-n+1)}.$$

On the other hand, evaluating the Mellin transform of $h$, we get
$$H(z)=\int_0^\infty h(t)t^{z-1}dt=\int_0^T x^{(n)}(T-t)t^{z-1}dt.$$
As $x^{(n)}$ is continuous and periodic, it follows that it is bounded, and we denote $\|x^{(n)}\|_\infty=\sup\limits_{t\in[0,T]}|x^{(n)}(t)|$. Hence,
$$|H(z)|\leq \|x^{(n)}\|_\infty \int_0^T|t^{z-1}|dt=\|x^{(n)}\|_\infty\int_0^Tt^{\Re(z)-1}dt$$
The integral from the right hand side of the inequality is convergent if and only if $\Re(z)>0$, and so the whole half-plane $\Re(z)>0$ is included in $S_h$.

Since $H(z)G(1-z)=0$, for any $z\in\mathbb{C}$ such that $z\in S_h$ and $1-z\in S_g$, it follows that
$$H(z)G(1-z)=0\qquad\forall z\in\mathbb{C}\textrm{ such that } n-\alpha<\Re(z)<1.$$
Since $G(1-z)\neq 0$ for any $z$ in the strip $n-\alpha<\Re(z)<1$, it follows that
$$H(z)=0 \qquad \forall z\in\mathbb{C}\textrm{ such that } n-\alpha<\Re(z)<1.$$
Applying the inverse Mellin transform, according to Proposition \ref{prop.Mellin.inversion} we obtain that the function $h$ is equal to $0$ almost everywhere on $(0,\infty)$, and therefore, $x^{(n)}$ is identically null (since it is continuous). Since all the derivatives of the function $x$, up to the order $n$, are $T$-periodic functions, it can be easily seen that the function $x$ is constant, which contradicts the initial assumptions. The proof of the theorem is now complete.
\qed\end{pf}

\begin{rem}
The proof of Theorem \ref{thm.main} relies on the definition of the Caputo derivative of order $\alpha\in(0,\infty)\setminus \mathbb{N}$ and the properties of the Mellin transform. When we talk about the usual integer-order derivatives, it is obvious that the derivative of a differentiable periodic function will be a periodic function of the same period. Therefore, Theorem \ref{thm.main} highlights a very important difference between fractional-order and integer-order derivatives.

We emphasize that the arguments from the proof of this theorem do not hold in the case when $\alpha\in \mathbb{N}$, i.e., when we deal with the usual integer-order derivatives. Indeed, if we assume that $\alpha\in \mathbb{N}$, then $n=[\alpha]+1=\alpha+1$ and the function $g$ defined in the proof of Theorem \ref{thm.main} is $g(t)=1$, for any $t\in(0,\infty)$. Since the integral $\int_0^\infty t^{z-1}dt$ is divergent, for any $z\in\mathbb{C}$, it follows that the Mellin transform of the function $g$ does not exist.
\end{rem}

\begin{rem}
In the papers \cite{Tavazoei-Haeri,Tavazoei}, the authors attempt to prove the same result as the one presented in Theorem \ref{thm.main}, but their arguments are flawed. One may clearly expect that at a certain point, the arguments provided in the proofs of the statements will fail for positive integer values of the order $\alpha$ (due to the facts described in the previous remarks). Nevertheless, a careful examination reveals that this is not the case in these papers. In fact, there is a mistake in the proof of Lemma 3 from \cite{Tavazoei-Haeri}, which constitutes one of the pre-requirements of the main results formulated in the Theorem from \cite{Tavazoei-Haeri} and Theorem 1 from \cite{Tavazoei}.

More precisely, the error comes from the fact that in the proof of the above mentioned Lemma 3, the areas $\bar{S}_i^+$ and $\bar{S}_i^-$ are wrongly considered to be independent of $p$. However, since $\bar{S}_i^+$ and $\bar{S}_i^-$ denote the area between the curve $(pT-t)^{\beta_i-1}\tilde{g}_i(t)$ and the $t$-axis, located above and bellow the axis, respectively (the interval $[0,T]$ is considered), they clearly depend on the values $p\in\mathbb{N}$. In fact, eq. (31) from \cite{Tavazoei-Haeri} should be written as
$$\int_0^T(pT-t)^{\beta_i-1}\tilde{g}_i(t)dt=\bar{S}_i^+(p)-\bar{S}_i^-(p).$$
 Therefore, based on (10) and (32) from \cite{Tavazoei-Haeri}, some simplifications lead us to:
$$p^{\alpha_i-\beta_i}\bar{S}_i^+(p)-(p-1)^{\alpha_i-\beta_i}\bar{S}_i^-(p)\leq 0 \leq (p-1)^{\alpha_i-\beta_i}\bar{S}_i^+(p)-p^{\alpha_i-\beta_i}\bar{S}_i^-(p),$$
where $\alpha_i<\beta_i$. This can be re-written as
$$\left(\frac{p-1}{p}\right)^{\beta_i-\alpha_i}\leq\frac{\bar{S}_i^+(p)}{\bar{S}_i^-(p)}\leq \left(\frac{p}{p-1}\right)^{\beta_i-\alpha_i}.$$
Passing to the limit when $p\rightarrow\infty$, we can simply obtain that
$$\lim_{p\rightarrow\infty}\frac{\bar{S}_i^+(p)}{\bar{S}_i^-(p)}=1$$
but it is impossible to conclude that $\bar{S}_i^+(p)=\bar{S}_i^-(p)$, for any $p\in\mathbb{N}$, which would lead us to the conclusion of Lemma 3.

In conclusion, it has to be underlined that the statements of Theorem 1 from \cite{Tavazoei} and the main Theorem from \cite{Tavazoei-Haeri} are correct, but the proofs of these statements rely on Lemma 3 from \cite{Tavazoei-Haeri}, which is incorrect.

We emphasize that the Mellin transform approach presented in this paper for the proof of Theorem \ref{thm.main} is essentially different from the one attempted in \cite{Tavazoei-Haeri,Tavazoei}.
\end{rem}

\begin{ex}
Considering $\alpha\in(0,1)$ and the function $x(t)=\sin(t)$, based on Theorem 6 from \cite{Ishteva} we obtain
$$^CD^\alpha \sin(t)= \frac{1}{2}t^{1-\alpha}\left[E_{1,2-\alpha}(it)+E_{1,2-\alpha}(-it)\right],$$
where $E_{\alpha,\beta}$ denotes the two-parameter Mittag-Leffler function defined by
$$E_{\alpha,\beta}(z)=\sum_{k=0}^\infty \frac{z^k}{\Gamma(\alpha k+\beta)}.$$
From the properties of the Mittag-Leffler function (see eq. (1.93) in \cite{Podlubny}) it follows that
$$^CD^\alpha \sin(t)= t^{1-\alpha}E_{2,2-\alpha}(-t^2).$$
For $\alpha\in(0,1)$, this function is not periodic (it can be verified by numerical simulations, see Fig. 1), even though it tends asymptotically to the periodic function $\sin\left(t+\frac{\alpha\pi}{2}\right)$. For $\alpha=1$, we have $ t^{1-\alpha}E_{2,2-\alpha}(-t^2)=E_{2,1}(-t^2)=\cos(t)$.

\begin{figure}[htbp]
\label{fig0}
\centering
\includegraphics*[width=0.7\linewidth]{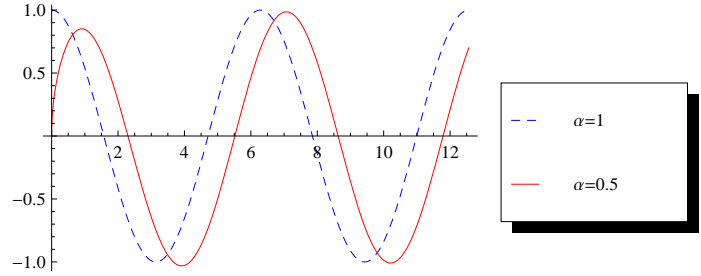}
\caption{The functions $\cos(t)$ (blue, dashed) and $t^{1-\alpha}E_{2,2-\alpha}(-t^2)$ (red) for $\alpha=0.5$.}
\end{figure}
\end{ex}

The following corollary (see also Theorem 2 from \cite{Tavazoei}) holds:

\begin{cor}\label{cor.RL.GL}
Let $\alpha\in (0,\infty)\setminus\mathbb{N}$ and $n=[\alpha]+1$. If $x:(0,\infty)\rightarrow\mathbb{R}$ is a non-constant $T$-periodic function of class $C^n$ on $(0,\infty)$, its Riemann-Liouville derivative $^{RL}D^\alpha x$ and its Grunwald-Letnikov derivative $^{GL}D^\alpha x$ cannot be $T$-periodic functions.
\end{cor}

\begin{pf}
Based on Proposition \ref{prop.relation}, we have
$$^{RL}D^\alpha x(t)=^{C}D^\alpha x(t) +\sum_{k=0}^{n-1}\frac{x^{(k)}(0^+)}{\Gamma(k-\alpha+1)}t^{k-\alpha},$$
where $n=[\alpha]+1$. Denoting $\|x^{(n)}\|_\infty=\sup\limits_{t\in[0,T]}|x^{(n)}(t)|$, we have:
$$
|^CD^{\alpha}x(t)|\leq\frac{1}{\Gamma(n-\alpha)}\int_0^t(t-s)^{n-\alpha-1}|x^{(n)}(s)|ds\leq \frac{\|x^{(n)}\|_\infty}{\Gamma(n-\alpha+1)}t^{n-\alpha},
$$
and it is easy to see that $^CD^{\alpha}x(0^+)=0$.

If there exists $k\in\{0,1,...,n-1\}$ such that $x^{(k)}(0^+)\neq 0$, it follows that
$$\lim_{t\rightarrow 0^+} {}^{RL}D^\alpha x(t)=\lim_{t\rightarrow 0^+}\sum_{k=0}^{n-1}\frac{x^{(k)}(0^+)}{\Gamma(k-\alpha+1)}t^{k-\alpha}=\pm\infty.$$
Hence, $^{RL}D^\alpha x(t)$ is unbounded in a neighborhood of $0$, and it cannot be periodic.

On the other hand, if $x^{(k)}(0^+)= 0$, for any $k\in\{0,1,...,n-1\}$, it follows that $^{RL}D^\alpha x(t)=^{C}D^\alpha x(t)$, and hence, based on Theorem \ref{thm.main}, it follows that $^{RL}D^\alpha x(t)$ cannot be $T$-periodic.

The conclusion about the Grunwald-Letnikov derivative $^{GL}D^\alpha x(t)$ can be easily derived from the equality $^{RL}D^\alpha x(t)=^{GL}D^\alpha x(t)$.
\qed\end{pf}

\section{Non-existence of periodic solutions in a class of fractional-order dynamical systems}

Let $\alpha=(\alpha_1,\alpha_2,...,\alpha_p)$, with $\alpha_i\in(0,\infty)\setminus\mathbb{N}$ for any $i\in\{1,2,...,p\}$ and $m\in\mathbb{N}$. We consider the following general class of fractional-order differential systems given in vector form:
\begin{equation}\label{sys.fractional}
D^{\alpha}\mathbf{x}(t)= \mathbf{f}(t,\mathbf{x}(t),\mathbf{x}'(t),...,\mathbf{x}^{(m)}(t)), \quad  t\geq 0,
\end{equation}
where $\mathbf{f} : [0,\infty)\times \mathbb{R}^{p(m+1)}\to \mathbb{R}^p$ is a given function and $$D^{\alpha}\mathbf{x}(t)=(D^{\alpha_1}x_1(t),D^{\alpha_2}x_2(t),...,D^{\alpha_p}x_p(t))^T$$ where $D^{\alpha_i}$, $i\in\{1,2,...,p\}$, is understood as one of the fractional-order derivatives (Riemann-Liouville, Caputo or Grunwald-Letnikov) given in Definition 1.

In the following, we denote $\mathbf{n}=([\alpha_1],[\alpha_2],...,[\alpha_p])\in\mathbb{N}^p$ and we say that the function $\mathbf{x}(t)=(x_1(t),x_2(t),...,x_p(t))^T$ is of class $C^{\mathbf{n}}(I,\mathbb{R}^p)$ if and only if $x_i(t)$ is of class $C^{n_i}(I,\mathbb{R})$, for any $i\in\{1,2,...,p\}$, where $I$ denotes a real interval.

The main results presented in the previous section lead us to the following consequence:
\begin{cor}[The non-autonomous case]\label{cor.nonautonomous}
Let $k\in\mathbb{N}^\star$ and assume that the function $\mathbf{f}$ is $T$-periodic with respect to its first argument. Then there are no non-constant $kT$-periodic solutions of class $C^{\mathbf{n}}$ of system (\ref{sys.fractional}).
\end{cor}

\begin{pf}
Assuming that there is a non-constant $kT$-periodic solution $\mathbf{x}(t)$ of class $C^{\mathbf{n}}$ of system (\ref{sys.fractional}), from the $T$-periodicity of the function $\mathbf{f}$ with respect to its first argument, it follows that $D^{\alpha}\mathbf{x}(t)$ is $kT$-periodic as well, which contradicts Theorem \ref{thm.main} or Corollary \ref{cor.RL.GL}.\qed
\end{pf}

In a very similar manner, the following result can be proved in the autonomous case:
\begin{cor}[The autonomous case]\label{cor.autonomous}
If the function $\mathbf{f}$ is constant with respect to its first argument (i.e., the system (\ref{sys.fractional}) is autonomous), then there are no non-constant periodic solutions of class $C^{\mathbf{n}}$ of system (\ref{sys.fractional}).
\end{cor}

\begin{rem}
Corollaries \ref{cor.nonautonomous} and \ref{cor.autonomous} hold even if some (but not all) of the differentiation orders $\alpha_1,\alpha_2,...,\alpha_p$ from system (\ref{sys.fractional}) are positive integers, i.e., some (but not all) of the fractional derivatives are replaced by integer-order derivatives.
\end{rem}

\begin{rem}
Even though, based on Corollary \ref{cor.autonomous}, exact periodic solutions do not exist in autonomous fractional-order systems, oscillatory behavior (limit cycles) has been observed by numerical simulations in many systems such as: a fractional-order Van der Pol system \cite{Barbosa}, fractional-order Chua and Chen's systems \cite{Cafagna-1,Cafagna-2}, a fractional-order R\"{o}ssler system \cite{Zhang-Zhou} and a fractional-order financial system \cite{Abd-Elouahab}.
\end{rem}

\begin{rem}
In the past few years, important existence results for solutions of periodic boundary value problems involving fractional differential equations have been obtained \cite{Belmekki-Nieto,Wei}, but the theoretical setting is essentially different from the one presented in Corollaries \ref{cor.nonautonomous} and \ref{cor.autonomous}.
\end{rem}

\begin{rem}
We emphasize that it is possible to obtain exact periodic solutions in impulsive fractional-order dynamical systems, by choosing the correct impulses at the right moments of time. For example, we may consider the impulsive dynamical system
\begin{align}\label{impulse}
\nonumber ^{C}D^{\alpha}\bold{x}(t)&= f(t,\bold{x}(t)), \qquad t>0, \;t\neq t_{k},  \\
\Delta \bold{x}\big|_{t=t_{k}}&= I_{k}(\bold{x}), ~~\quad\qquad k\in\mathbb{Z}^\star_+\\
\nonumber \bold{x}(0)&= \bold{x}_{0},
\end{align}
where $\alpha\in(0,1)$, $f:[0,\infty)\times \mathbb{R}^n\to \mathbb{R}^n$ is a given function, $\bold{x}_{0}\in\mathbb{R}^n$, the sequence of times $(t_k)_{k\in\mathbb{Z}_+}$ is strictly increasing, $t_0=0$, $\Delta \bold{x}|_{t=t_{k}}=\bold{x}(t_{k}^{+})-\bold{x}(t_{k}^{-})$, where
$\bold{x}(t_{k}^{+})$ and $\bold{x}(t_{k}^{-})$ represent the right and left limits of $\bold{x}(t)$ at $t=t_{k}$, and $I_k$ are the impulsive operators. We consider that the following assumptions are satisfied:
\begin{itemize}
  \item The function $f(t,\cdot)$ is $T$-periodic.
  \item Let $p\in\mathbb{Z}_+$ such that $[0,T]\cap (t_k)_{k\in\mathbb{Z}_+}=\{t_0,t_1,t_2,...,t_p\}$. We assume: $t_{k+p}=t_k+T$, for any $k\in\mathbb{Z}_+$.
\end{itemize}
Defining the impulsive operators by
$$I_{k}(\bold{x})=-\frac{1}{\Gamma(\alpha)}\int_{t_{k-1}}^{t_{k}}(t_{k}-s)^{\alpha-1}f(s,\bold{x}(s))ds,\qquad\forall k\in\mathbb{Z}_+^\star$$
ensures that the solution of (\ref{impulse}) is $T$-periodic. Theoretical details and several applications will be given in a future paper.

\end{rem}

\section{An application to fractional-order neural networks}

The main advantage of fractional-order mathematical models in comparison with classical integer-order models is that fractional derivatives provide an excellent tool for the description of memory and hereditary properties of various processes. In fact, fractional-order systems have infinite memory. Taking into account these facts, it is easy to see that the incorporation of a memory term (in the form of a fractional derivative or integral) into a neural network model is an extremely important improvement.

Based on this idea, the common capacitor from the continuous-time integer-order Hopfield neural network is replaced by a generalized capacitor, called fractance, giving birth to the so-called fractional-order Hopfield neural network model.

The fractional-order formulation of artificial neural network models is also justified by research results concerning biological neurons. In the recent paper \cite{Lundstrom} it has been pointed out that fractional differentiation provides neurons with a fundamental and general computation ability that can contribute to efficient information processing, stimulus anticipation and frequency-independent phase shifts of oscillatory neuronal firing, emphasizing once again the utility of developing and studying fractional-order mathematical models of neural network dynamics.

The analysis of fractional-order artificial neural networks is a very recent and promising research topic. The first papers discussing fractional-order neural network models \cite{Arena,Matsuzaki,Petras,Boroomand} report on results of numerical simulations, especially on the numerical evidence of limit cycles and chaotic phenomena.

We consider a simple fractional-order neural network of two neurons, given by the following autonomous system of fractional-order differential equations
\begin{equation}\label{ex.2d}
\left\{\begin{array}{l}
^CD^{\alpha} x_1(t)=-x_1(t)+2 \tanh(x_1(t))-0.5\tanh(x_2(t))\\
^CD^{\alpha} x_2(t)=-x_2(t)+\tanh(x_1(t))+2 \tanh(x_2(t))
\end{array}
\right.
\end{equation}
where the fractional order is $\alpha=0.5$.

\begin{figure}[htbp]
\label{fig2}
\centering
\includegraphics*[width=0.6\linewidth]{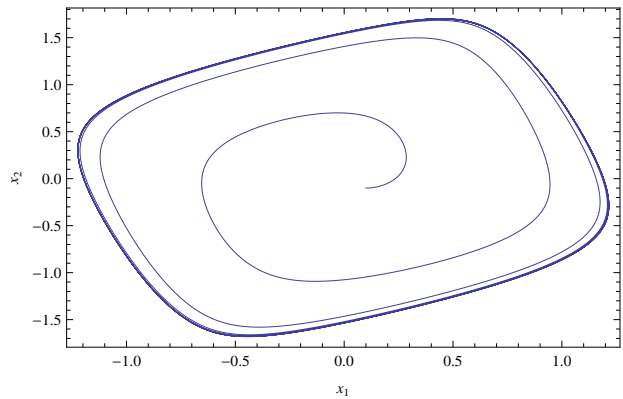}
\caption{The trajectory of system (\ref{ex.2d}), with the initial condition $(x_0,y_0)=(0.1,0.1)$ converges to an asymptotically stable limit cycle.}
\end{figure}

\begin{figure}[htbp]
\label{fig2}
\centering
\includegraphics*[width=0.8\textwidth]{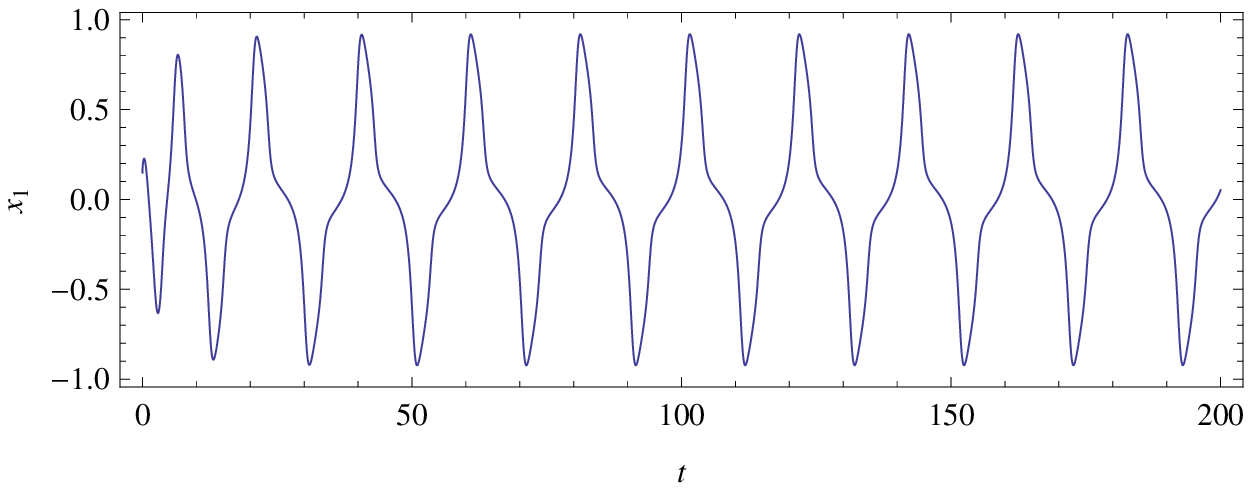}\\
\includegraphics*[width=0.8\textwidth]{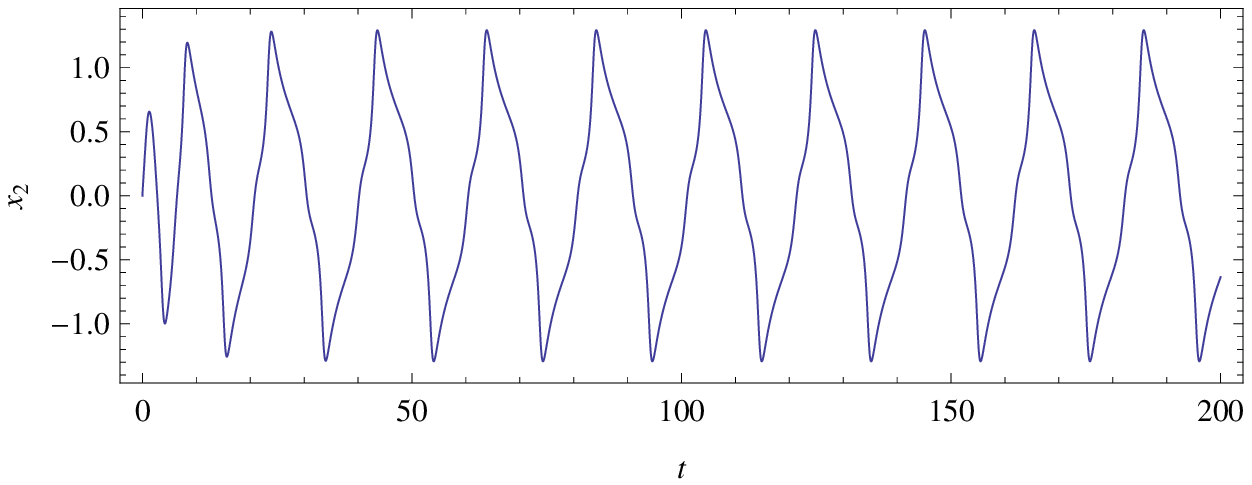}
\caption{The solution $(x_1(t),x_2(t))$ of system (\ref{ex.2d}) versus time, with the initial condition $(x_0,y_0)=(0.1,0.1)$.}
\end{figure}

Numerical simulations suggest that there exists an attracting limit cycle for system (\ref{ex.2d}) (see Figures 2 and 3). However, based on Corollary \ref{cor.autonomous}, the attracting limit cycle is not an exact periodic solution of system (\ref{ex.2d}). We refer to \cite{ES-IJCNN-2011}
 for additional results concerning system (\ref{ex.2d}).

\section{Conclusions and directions for future research}

In this paper, it has been shown that the fractional-order derivative of a periodic function cannot be a periodic function of the same period, underlining a very important difference between integer and fractional-order derivatives and explaining the absence of periodic solutions in a class of fractional-order dynamical systems. Since this remarkable property holds for three widely used types of fractional derivatives (Riemann-Liouville, Caputo and Grunwald-Letnikov), it may be interesting to investigate if it is also true for other known definitions of fractional derivatives, or if it is possible to give a new formal definition of a fractional-order derivative which preserves periodicity. For example, in section 3.1 of \cite{Butzer-Westphal}, a generalized definition was presented for fractional derivatives of periodic functions, using Fourier expansions, and several results concerning periodicity were also given. Furthermore, two kinds of Weyl fractional derivatives defined for periodic functions \cite{Hilfer}, namely the  Weyl-Liouville and Weyl-Marchaud fractional derivatives, are also worth mentioning.



\end{document}